\documentclass[final,leqno,onefignum,onetabnum]{siamltex1213}

\usepackage{graphicx}
\usepackage{amssymb}
\usepackage{amsmath}

\def\xR{\mathbb{R}}

\def\xi{\mathbf{i}}

\def\xn{\mathbf{n}}
\def\xx{\mathbf{x}}
\def\xv{\mathbf{v}}
\def\dd{\textnormal{d}}

\def\crit{\textnormal{crit}}

\newtheorem{Th}{Theorem}
\newtheorem{Lemma}{Lemma}[section]

\newtheorem{Rem}{Remark}[section]

\title{OOID GROWTH: UNIQUENESS OF STEADY STATE, SMOOTH SHAPES IN 2D} 

\author{Andr\'as A. Sipos \thanks{MTA-BME Morphodynamics Research Group, Budapest, Hungary \& Department of Mechanics, Materials and Structures, Budapest University of Technology and Economics, Budapest, Hungary.
(\email{siposa@eik.bme.hu}).}}

\begin{document}
\maketitle
\slugger{mms}{xxxx}{xx}{x}{x--x}

\begin{abstract}
Evolution of planar curves under a nonlocal geometric equation is investigated. It models the simultaneous contraction and growth of carbonate particles called ooids in geosciences. Using classical ODE results and a bijective mapping we demonstrate that the steady parameters associated with the physical environment determine a unique, time-invariant, compact shape among smooth, convex curves embedded in $\xR^2$. It is also revealed that any time-invariant solution possesses $D_2$ symmetry. The model predictions remarkably agree with ooid shapes observed in nature. 
\end{abstract}



\pagestyle{myheadings}
\thispagestyle{plain}

\section{Introduction}
A geometric, non-local PDE is considered to model the shape evolution of mm-sized carbonate particles called \emph{ooids}. They form in shallow tropical coastal waters and are widely investigated as important markers of coastal environments in the geological past. In \cite{Sipos2018} a simple, two-dimensional model of ooid growth is introduced as a natural extension of the global model of \cite{Trower2017}. The latter hypothesized and experimentally verified that ooids reflect a precious balance between increase and reduction of the grain's net volume. In the pointwise model of \cite{Sipos2018} the curve representing the shape is mapped in the local normal direction. The speed of the motion is driven by three, well-distinguished physical processes. These are \emph{chemical precipitation} leading to radial accumulation of material, \emph{abrasion} of the grain due to collisions with the seabed and sliding \emph{friction}, which takes effect at shallow shores. \cite{Sipos2018}
presents numerical evidence about time-invariant solutions, and remarkable resemblance to cross-sections of real ooids is found. Furthermore, a hypothesis about the bold intermediate layers widely observed in ooid cross sections is established. The present paper is devoted to the rigorous investigation of the existence and uniqueness of time-invariant solutions of the model introduced in \cite{Sipos2018}.

Generally speaking, shape evolution of particles is widely investigated both in the mathematical and in the geoscientific literature (e.g. \cite{Bloore1977,Domokos2012} and the citations therein). Most of the treated models are local ones, i.e. the evolution is determined by some pointwise law; the \emph{curve-shortening flow} \cite{Grayson1989} is a good example for such a model in two spatial dimensions. Perhaps investigation of ancient solutions under some prescribed flow (e.g. \cite{Hamilton2010}) is the closest to our problem, however, here the simultaneous presence of growth and reduction of the shape makes it straightforward to seek compact, time-invariant shapes under the flow.

\section{The model and the main results}
\label{sec:main}

Shape evolution might be interpreted as a process that moves any point $\xx=(x,y)$ of a closed, non-self-intersecting curve $\Gamma$ embedded in $\xR^2$ to the inward normal direction $\xn$ with a speed $\xv$ that depends on intrinsic features of the curve and parameters characterizing the physical environment. The geometric evolution equation to model ooid growth introduced in \cite{Sipos2018} reads
\begin{equation}
\label{eq:geom_PDE}
\xx_t=\xv=c_3\left(-1+c_1A\kappa +c_2Ay\cos\gamma\right)\xn,
\end{equation} 

\noindent where $A$ is the - time dependent - area enclosed by $\Gamma$ and the $t$ subscript refers to differentiation with respect to time. $\kappa$ and $\gamma$ stand for the \emph{curvature} and the \emph{turning angle}, respectively (cf. Figure 1). Any parametrization of $\Gamma$ makes the quantities $\kappa$ and/or $\gamma$ in (\ref{eq:geom_PDE}) to be dependent on derivatives with respect to the parametrization, which reveals that (\ref{eq:geom_PDE}) is in fact a parabolic PDE. 

Following the lead of \cite{Sipos2018} we assume that $\Gamma$ possesses a unique maximal diameter (line $e$ between points P and P' in Fig. \ref{Fig:01}.), which is designated to be the $x$ axis of an orthonormal basis located at the middle point of the PP' segment. This assumption makes friction well-defined in the model: the affine law in the third term of (\ref{eq:geom_PDE}) produces abrasion in points far away from axis $x$, which is expected from a sliding motion parallel to it. Nonetheless, other formulations of sliding friction might be physically pausible. However -- as it is pointed out in \cite{Sipos2018} -- the two, area-dependent terms in (\ref{eq:geom_PDE}) result in a second-order approximation of the change in the net volume which justifies our choice, both in formulating the frictional law and making the collisional and frictional terms proportional to the enclosed area (a.k.a. the net volume in 2D).

Now $\gamma$ denotes the angle between the $x$ direction and the local tangent to the curve. $c_1$, $c_2$ and $c_3$ are positive real parameters associated with the physical environment and they are assumed to be time-independent during the course of shape evolution. Their dimensions are length$^{-1}$, length$^{-3}$ and lentgh/time, respectively. The three key physical processes driving the evolution can be easily identified: in the brackets the first, negative term stands for \emph{growth}, in the second term \emph{abrasion} is assumed to be a curvature-driven process and finally the affine term is associated with \emph{friction}. As abrasion and friction are proportional to mass (and growth is not), the last two terms in 2D depends on the global quantity $A$. As the first term is negative and the the other two are positive, we seek compact invariant shapes (denoted to $\Gamma^*$) that fulfill
\begin{equation}
\label{eq:steady_state}
-1+c_1A\kappa +c_2Ay\cos\gamma=0, \qquad\qquad \forall\xx\in\Gamma^*.
\end{equation} 

\noindent Note that $\Gamma^*$ is independent of $c_3$ as it scales solely the time and cannot be reconstructed by pure observation of the shape. In \cite{Sipos2018} it is demonstrated, that although the friction term contains orthogonal affinity, ellipses are not invariant solutions. In this paper we show that among smooth, convex curves any time-invariant shape under the above-defined flow must possess $D_2$ symmetry (i.e. its symmetry group is generated by a non-square rectangle). Furthermore, for a given parameter pair $(c_1,c_2)$ the invariant shape is unique.

\begin{Th}
\label{thm1}
Let the parameters in (\ref{eq:steady_state}) be time-invariant, positive constants ($c_1 >0$ and $c_2\geq 0$). Then any smooth, convex, time-invariant curve $\Gamma^*$ satisfying (\ref{eq:steady_state}) in all of its points and embedded in $\xR^2$ possesses $D_2$ symmetry.
\end{Th}

\begin{Th}
\label{thm2}
The smooth, convex, time-invariant curves under the flow in (\ref{eq:steady_state}) are uniquely determined by $c_1$ and $c_2$, and for any positive values of these parameters there exists a $\Gamma^*$ curve.
\end{Th}

We prove the first theorem in Section \ref{sec:Prop1}, where we assume that $A$ is known \emph{a-priori}, this case is refereed to as \emph{local equation} to distinguish it from the general, \emph{non-local equation}. Section \ref{sec:Prop2} is devoted for the proof of Theorem 2. Finally, conclusions are drawn.

\section{The local equation}
\label{sec:Prop1}

For a moment let us assume that the area of the invariant curve is known \emph{a-priori}. (This assumption can be justified by imagining the flow with fixed parameters to be run until a steady state. If it happens, the area can be measured.) Without loss of generality we consider solely the curve segment $\bar\Gamma$ between the leftmost point P and the one that possesses a horizontal tangent and a positive $y$ coordinate (point Q). In order to simplify the derivations we use several parametrizations of the curve segment in the sequel: parametrization with respect to the arc length (natural parametrization), to the $y$ coordinate and to the $\gamma$ turning angle, respectively. Derivatives with respect to the parametrization is denoted by lower indexes.

\begin{figure}[!ht]
\centering
\includegraphics[width=0.85\textwidth]{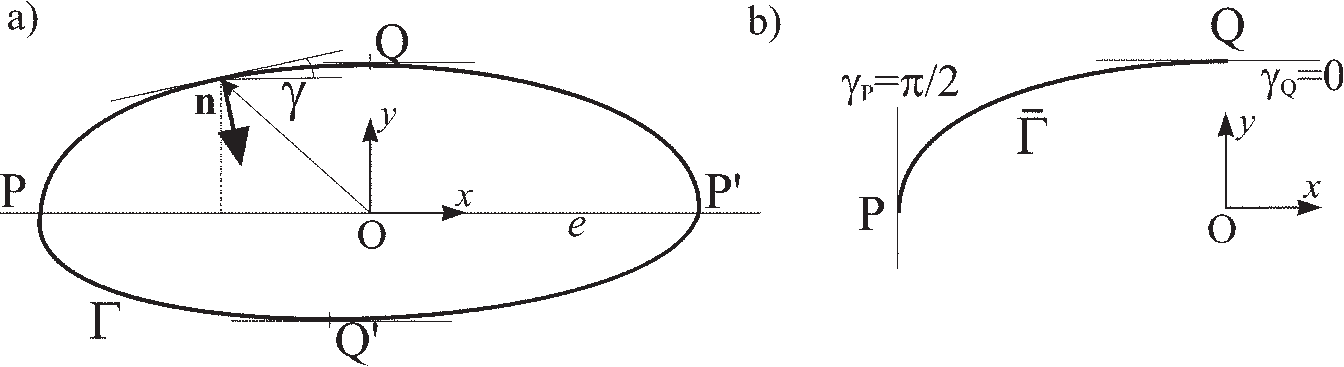}
\caption{Notations}
\label{Fig:01}
\end{figure}  

In the local equation $A$ is fixed, hence it is convenient to introduce $\hat{c}_1=c_1A$ and $\hat{c}_2=c_2A$, which renders (\ref{eq:steady_state}) into
\begin{equation}
\label{eq:local_02}
-1+\hat{c}_1\kappa+\hat{c}_2y\cos\gamma=0.
\end{equation}

\begin{Lemma}
\label{lemma3}
For fixed parameters $\hat{c}_1$ and $\hat{c}_2$ there is at most one curve segment $\bar\Gamma$ that fulfills equation (\ref{eq:local_02}) in all of its interior points.
\end{Lemma}
\begin{proof}
For a moment we reconsider the natural parametrization of the curve with the arch length parameter $s$. If $s$ increases clockwise, then $x_s=\cos\gamma$ and $y_s=\sin\gamma$. Recall that the derivative of the slope respect to the arch length equals the curvature. The chain rule yields
\begin{eqnarray}
\label{eq:DE_kappa_gamma}
\kappa(y):=-\gamma_s=-\gamma_y y_s=-\gamma_y\sin\gamma=(\cos\gamma)_y,
\end{eqnarray}

\noindent where the negative sign indicates, that $\gamma(y)$ is decreasing between points P and Q (Fig. \ref{Fig:01}. b)). For brevity let 
\begin{align}
\label{eq:def_q}
q:&=\frac{\hat{c}_2}{2\hat{c}_1}.
\end{align}

\noindent Introducing $\phi(y):=\hat{c}_1\cos(\gamma(y))$, substituting (\ref{eq:DE_kappa_gamma}) into (\ref{eq:local_02}) and applying (\ref{eq:def_q}) yield
\begin{equation}
\label{eq:local_03}
-1-\hat{c}_1\gamma_y\sin\gamma+\hat{c}_2y\cos\gamma=-1+\phi_y+2qy\phi=0,
\end{equation}

\noindent which is a first order, linear ODE. Classical results of ODE theory provide existence and uniqueness for $\phi(y)$ and consequently for $\kappa=\hat{c}_1^{-1}\phi_y$. In specific, the integrating factor $u(y)$ reads
\begin{align}
u(y):=\exp\int 2qy\dd y=\exp qy^2.
\end{align}

Recall that in our model $\gamma(0)=\pi/2$ which yields $\phi(0)=0$. Hence, the Cauchy problem in (\ref{eq:local_03}) with the initial condition $\phi(0)=0$ possesses the following solution:
\begin{align}
\label{eq:def_phi}
\phi(y)=\frac{\int_0^yu(\eta)\dd\eta}{u(y)}=\int_0^y\exp\left(q(\eta^2-y^2)\right)\dd\eta.
\end{align}

The curvature function readily follows from (\ref{eq:DE_kappa_gamma})
\begin{equation}
\label{eq:kappa(delta)}
\kappa(y)=\frac{1}{\hat c_1}\phi_y=\frac{1}{\hat c_1}\left(1-2qy\int_0^y\exp\left(q(\eta^2-y^2)\right)\dd\eta\right).
\end{equation}

Detailed investigation of the properties of $\kappa(y)$ is needed for further development. From the r.h.s. of (\ref{eq:local_03}) the first, second and third derivatives of $\phi(y)$ are obtained
\begin{align}
\label{eq:phi_y}
\phi_y&=-2qy\phi+1,\\
\label{eq:phi_yy}
\phi_{yy}&=-2qy\phi_y-2q\phi,\\
\label{eq:phi_yyy}
\phi_{yyy}&=-2qy\phi_{yy}-4q\phi_y.
\end{align}

Using these, the following properties of $\kappa(y)$ can be settled: 
\begin{enumerate}
	\item $\phi(y)$ and consequently $\kappa(y)$ are $C^\infty$, which is obvious from (\ref{eq:def_phi}) and (\ref{eq:kappa(delta)}).
	\item $\kappa(0)$ is positive and equals $\hat{c}_1^{-1}$. As $\kappa(0)=\hat c_1^{-1}\phi_y(0)$ the claim follows.
	\item $\kappa(y)$ has a local maximum at $y=0$. Firstly, $\kappa_{y}(0)=\hat c_1^{-1}\phi_{yy}(0)=-2\hat c_1^{-1}q\phi(0)=0$ indicates that $y=0$ is indeed a critical point. Secondly, $\kappa_{yy}(0)=\hat c_1^{-1}\phi_{yyy}(0)=-\hat c_1^{-1}4q\phi_y(0)=-4\hat c_1^{-1}q <0$, which demonstrates, that the critical point at $y=0$ is a maximum.
		
	\item $\kappa(y)\rightarrow 0$ as $y\rightarrow\infty$. Using l'Hopital's rule we have
\begin{align}
\nonumber\lim_{y\rightarrow\infty}y\phi(y)=\lim_{y\rightarrow\infty}\frac{y\int_0^y\exp(q\eta^2)\dd\eta}{\exp(qy^2)}=\lim_{y\rightarrow\infty}\frac{(1+2qy^2)\exp(qy^2)}{(2q+4q^2y^2)\exp(qy^2)}=\frac{1}{2q}.
\end{align}

\noindent This result, eq. (\ref{eq:phi_y}), and the fact that $\hat{c}_1$ and $q$ are fixed parameters yield the desired result as  
\begin{align}
\nonumber\lim_{y\rightarrow\infty}\kappa(y)=\lim_{y\rightarrow\infty}\frac{1}{\hat{c}_1}\left(-2qy\phi+1\right)=0.
\end{align}
		
	\item There is exactly one point, denoted to $y_0$, where $\kappa(y)$ vanishes and $y_0$ solely depends on $q$. By definition $\phi(y)\geq 0$ with $\phi(0)=0$. An analogous argument to point (4) above shows that $\lim_{y\rightarrow\infty}\phi(y)=0$. The positivity of $\phi(y)$ and $q$ yield that $\phi_{yy}$ in (\ref{eq:phi_yy}) is negative at any critical point for $0<y<\infty$. Hence, as $\phi(y)$ is smooth, it follows, that there is one, and only one point, at which $\phi_y=0$, hence $y_0$, at which $\kappa(y)$ vanishes exists and it is unique.
	
	\item There is no local extrema for $\kappa(y)$ between $0<y<y_0$, thus it is monotonic in this range. Note that points (2) and (5) yield $\kappa(y)>0$ and subsequently $\phi_y>0$ for $0<y<y_0$. Now, as $\phi(y)$ is positive, (\ref{eq:phi_yy}) shows that $\phi_{yy}<0$ in $0<y<y_0$. Hence, $\kappa_y=\hat{c}_1^{-1}\phi_{yy}$ is strictly negative excluding any critical point between $y=0$ and $y=y_0$.

\end{enumerate}

To realize an invariant shape $\Gamma^*$ we need $\gamma(y)$ itself. By the virtue of eq. (\ref{eq:DE_kappa_gamma})
\begin{equation}
\label{eq:gamma_from_kappa}
\gamma(y)=\arccos\left(\int_0^y\kappa(\eta)\dd\eta\right).
\end{equation}

\noindent Since $\arccos(.)$ is monotonic decreasing in $[0,1]$, the area below the solution function $\kappa(y)$ determines a unique $\gamma(y)$. In other words $\hat{c}_1$ and $q$ (or $\hat{c}_1$ and $\hat{c}_2$) determine a unique steady state curve for eq. (\ref{eq:local_02}); We aim to determine the parameter range, where the curve is smooth. Apparently, if the area under $\kappa(y)$ between $0\leq y\leq y_0$ exceeds 1, then we can construct a smooth shape: at the unique $\bar y<y_0$ the area below $\kappa(y)$ equals 1, i.e. this corresponds to point Q with a tangent parallel to the axis $x$. This solvability condition can be derived explicitly as follows. For a smooth shape we need
\begin{align}
\label{eq:solvability_1}
\int_0^{y_0}\kappa(\eta)\dd\eta=\frac{1}{\hat{c}_1}\phi(y_0)\geq 1.
\end{align}

Let us introduce $\zeta:=\sqrt q\eta$ and $z:=\sqrt q y$. After changing variables (\ref{eq:solvability_1}) reads
\begin{align}
\label{eq:solvability_2}
\frac{\sqrt{2}}{\sqrt{\hat{c}_1\hat{c}_2}}\int_0^{z}\exp(\zeta^2-z^2)\dd\zeta\geq 1,
\end{align}

\noindent hence the fixed parameters are required to fulfill
\begin{align}
\label{eq:solvability_3}
\sqrt{\hat{c}_1\hat{c}_2}\leq\sqrt{2}\max_{z\geq 0}\int_0^{z}\exp(\zeta^2-z^2)=:\Psi,
\end{align}

\noindent where the upper bound $\Psi$ is approximated numerically as $\Psi\approx 0.765$. If the condition in (\ref{eq:solvability_3}) is met, then from (\ref{eq:solvability_1}) $\phi(\bar{y})=\hat{c}_1$ follows.

For $0\leq y\leq \bar y$ the connection between $\gamma$ and $y$ is one to one, thus we can draw the \emph{physical realization}. For cases, at which condition (\ref{eq:solvability_1}) fails, the physical shapes are non-smooth (in fact, they become concave as the curvature flips sign above $y_0$ and there is no other zero for $\kappa(y)$). As we have seen, $\kappa(0)$ depends solely on $\hat{c}_1$ and for fixed $q$ the value of $y_0$ is fixed, too. The definition in (\ref{eq:def_q}) and the solvability condition in (\ref{eq:solvability_3}) lead to the conclusion that for any fixed $q$ there exists a $\hat{c}_{1,\crit}:=\Psi\sqrt{2q}$. For $0<\hat{c}_1\leq\hat{c}_{1,\crit}$ the invariant, smooth curve segment $\bar\Gamma^*$ is unique otherwise, there is no such solution.
\end{proof}

For further convenience at a fixed value of $q$ we introduce the set
\begin{equation}
\label{def:chi_q}
\chi_q:\left\{\hat{c}_1 \quad| \quad 0<\hat{c}_1\leq\hat{c}_{1,\crit} \right\}.
\end{equation}

Note that to have a nonzero measure of $\chi_q$ one needs $\hat{c}_1>0$. It follows, that for $\hat{c}_1>\hat{c}_{1,\crit}$ the integral on the left-hand-side of (\ref{eq:solvability_1}) is smaller than one, which means the associated curve cannot have a horizontal tangent at any point. Having assumed convex, smooth curves this parameter-range is not in our interest. In case $\hat{c}_1 \in \chi_q$ the shape can be realized.

\begin{Lemma}
\label{lemma4}
The closed, non-intersecting curve $\Gamma^*$ obtained by reflections of the curve segment $\bar\Gamma^*$ with respect to the axes $x$ and $y$ is a $C^\infty$ curve.
\end{Lemma}

\begin{proof}
By smoothness of the $\kappa(y)$ function $\bar\Gamma^*$ is smooth in its interior points. We need to prove, that $\Gamma^*$ is $C^\infty$ in points $P,P',Q$ and $Q'$. Without loss of generality, we show smoothness at points $P$ and $Q$, it follows by symmetry for $P'$ and $Q'$. For point $P$ observe that from (\ref{eq:def_phi}) $\phi(y)=-\phi(-y)$ follows, i.e. $\phi(y)$ is an odd function. Following (\ref{eq:phi_y})-(\ref{eq:phi_yyy}) it is straightforward to show, that derivatives of $\phi$ with respect to $y$ at $y=0$ fulfills:
\begin{itemize}
	\item odd derivatives are finite,
	\item even derivatives vanish.
\end{itemize}
We conclude, that $\phi(y)$ is analytic at $y=0$ which demonstrate that the curve is smooth at point $P$.
 
For point $Q$ re-parametrization of $\Gamma^*$ is essential as a parametrization with respect to $y$ is not one-to-one for $\Gamma^*$. Let the curve be parametrized with respect to its arch length $s$ such way, that at point $Q$ there is $s=0$ and $s$ increasing clockwise. By reflection we have $y(s)=y(-s)$. Let $\tilde\phi(s)$ denote the extension of $\phi(y)$ after the re-parametrization of $\Gamma^*$. In specific, after reflection we find
\begin{align}
\tilde\phi(s)=\phi(y(s))=\phi(y(-s))=\tilde\phi(-s),
\end{align}

\noindent which shows, that $\tilde\phi$ is an even function at point $Q$. In specific, at point $Q$ we have $\tilde\phi(0)=\phi(\bar{y})=\hat{c}_1$. Employing (\ref{eq:phi_y})-(\ref{eq:phi_yyy}) and the chain rule we find, that derivatives of $\tilde\phi$ with respect to $s$ at $s=0$ follow the pattern:
\begin{itemize}
	\item odd derivatives vanish,
	\item even derivatives are finite.
\end{itemize}
Once again, we conclude, that $\tilde\phi(s)$ is analytic at $s=0$, hence the smoothness of the curve at point $Q$ follows.
\end{proof}


\begin{proof}[Proof of Theorem \ref{thm1}]
By Lemma 3.1 positive parameters $\hat{c}_1$ and $\hat{c}_2$ determine a unique curve segment $\bar\Gamma^*$ with vertical tangent at point $P$ and horizontal tangent at point $Q$ iff $0<\hat{c}_1\leq \hat{c}_{1,\crit}$. By Lemma 4 reflections of $\bar\Gamma^*$ with respect to $x$ and $y$ produce a closed, convex, smooth curve. Finally, assumption of a single maximal diameter $e$ for $\Gamma^*$ implies uniqueness.

Solution of the local equation establish a solution for the non-local case (eq. \ref{eq:steady_state}), too. To see this, let us fix the two parameters, $\hat{c}_1$ and $\hat{c}_2$, follow the lines in this section to obtain a steady state solution $\Gamma^*$. In case there exists such a solution, measure area $A$ enclosed by the curve. It simply delivers the parameters of the non-local equation via $c_1=\hat{c}_1/A$ and $c_2=\hat{c}_2/A$. In the other way round, if one knows a time-invariant solution of the non-local equation, calculation of the parameters in the local is straightforward. These observations imply that a smooth solution of the non-local case must possess $D_2$ symmetry, too. 
\end{proof}

\begin{Rem}
For $\hat{c}_2=0$ we have $q=0$ and hence $\kappa(y)\equiv 1/(\hat{c}_1)=\kappa(0)$ implying that in this case the time-invariant shape is a circle. As the term of friction (the one with parameter $\hat{c}_2$) represents an affine flow, in the general case (i.e. $\hat{c}_2\neq 0$) $\kappa(0)$ represents the maximal curvature of the curve.
\end{Rem}

In the next Section we investigate the connection between the local and non-local models via the relations between their parameters.

\section{The non-local equation}
\label{sec:Prop2}
We turn to investigate steady state solutions of the non-local equation (\ref{eq:steady_state}). As we found that the symmetry group of any invariant solution is $D_2$, we keep investigating a curve segment $\bar\Gamma$ (c.f. Figure 1.). To investigate uniqueness of solutions in (\ref{eq:steady_state}) let us assign $(\hat{c}_1,\hat{c}_2)$ and $(c_1,c_2)$ if they result in and identical time-invariant curve of the proper model. In this sense we can talk about a \emph {mapping between the parameter spaces}. 

Observe that parameter $q$ in eq. (\ref{eq:def_q}) is invariant under this map because
\begin{align}
\frac{\hat{c}_1}{\hat{c}_2}=\frac{Ac_1}{Ac_2}=\frac{c_1}{c_2}.
\end{align}
In order to facilitate this observation, instead of $\hat{c}_2$ and $c_2$ we use $q$ as one of the parameters in the problem. Based on (\ref{eq:solvability_3}) and (\ref{def:chi_q}), in the local model only $\hat{c}_1\in\chi_q$ produces a smooth curve. For a fixed value of $q$ let the map $M$ defined as
\begin{eqnarray}
\label{eq:1D_map}
M:\chi_q \rightarrow \xR^+ \qquad\textrm{with}\qquad \hat{c_1} \mapsto c_1.
\end{eqnarray}

Our program is to show that $M$ is injective and surjective, thus it is bijective implying smooth solution curves of the non-local equation are unique as we had uniqueness of solutions for eq. (\ref{eq:local_03}).
\begin{Lemma}
Map $M$ is injective.
\end{Lemma}
\begin{proof}
As we have seen, $\hat{c}_1\in\chi_q$ results in a smooth curve enclosing some positive area $A$. Based on our construction, $c_1(\hat{c}_1):=\hat{c}_1A^{-1}$ can be readily computed. It means, injectivity of $M$ follows from the strict monotonicity of the $c_1(\hat{c}_1)$ function over $\chi_q$. To prove this, let us consider two smooth solutions (at a fixed value of $q$) of the local equation in (\ref{eq:local_03}) identified by the letters $i$ and $j$. Their parameters are related as
\begin{equation}
\label{eq:c_connection}
\hat{c}_1^j=(1+\varepsilon)\hat{c}_1^i, 
\end{equation}

\noindent where without loss of generality $\varepsilon>0$. By the virtue of eq. (\ref{eq:kappa(delta)}) is is clear, that not only the parameters, but the $\kappa(y)$ functions of the time-invariant curve segments $\bar\Gamma_i^*$ and $\bar\Gamma_j^*$ fulfill
\begin{equation}
\label{eq:kappa_connection}
\kappa^j(y)=\frac{1}{1+\varepsilon}\kappa^i(y) \qquad \forall y:\int_0^y\kappa{\eta}\dd\eta\leq 1.
\end{equation}

\begin{figure}[!ht]
\centering
\includegraphics[width=0.95\textwidth]{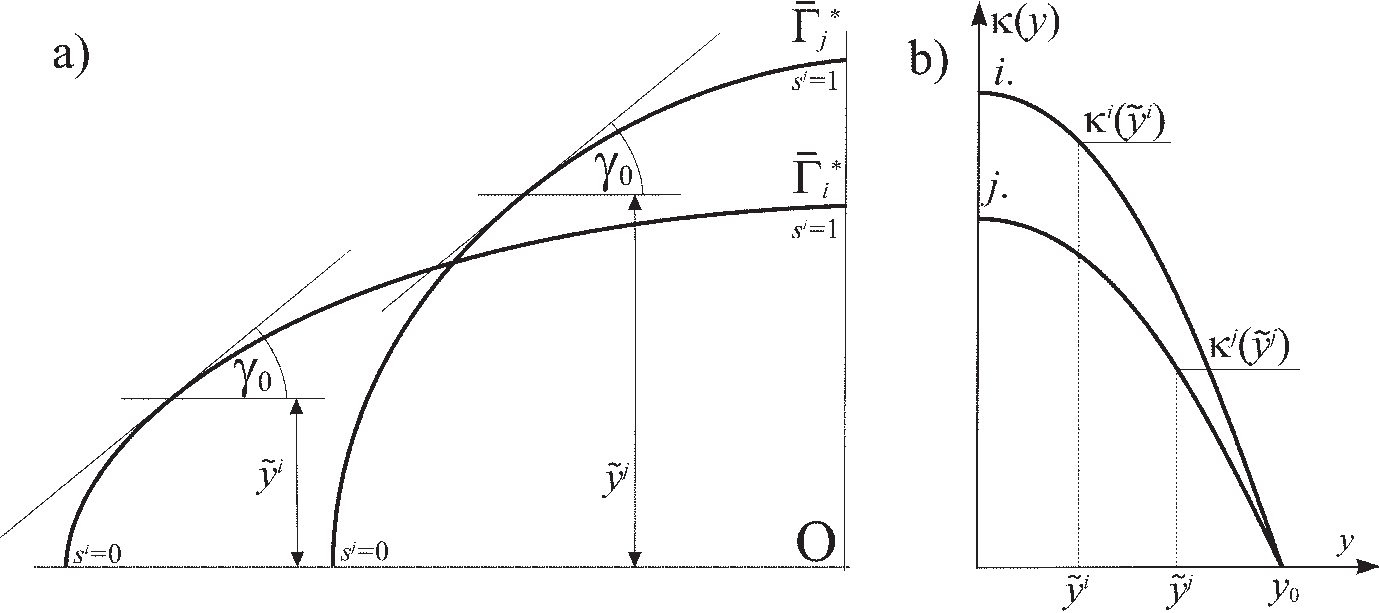}
\caption{Comparison of two, steady state curve-segments, $\bar\Gamma_i^*$ and $\bar\Gamma_j^*$ at fixed $q$. a) depicts the two segments and denote the arbitrary point-pair with a fix $\gamma_0$, which is used to determine the relation between the areas under the curve-segments. b) the graphs of $\kappa(y)$ curvature functions for $\bar\Gamma_i^*$ and $\bar\Gamma_j^*$, respectively.}
\label{Fig:02}
\end{figure}  

We choose two points along $\bar\Gamma_i^*$ and $\bar\Gamma_j^*$ (Fig. \ref{Fig:02}.), one for each, such way that their turning angles are identical. The common angle is denoted to $\gamma_0$ and the $\tilde{(.)}$ sign refers to any quantity evaluated at these points (e.g. $\tilde{y}^i$ is the parameter of curve $i$ at the chosen point along $\bar\Gamma_i^*$). As $\gamma(y)$ is monotonic along $\bar\Gamma$, the position of the two points is well-defined. As it is demonstrated in Section 3, $\kappa(y)$ and $\gamma(y)$ are related via (\ref{eq:gamma_from_kappa}), thus for our two curves we find that
\begin{equation}
\label{eq:gamma_connection}
\int_0^{\tilde{y}^i}\kappa^i(\eta)\dd\eta=\cos(\gamma_0)=\int_0^{\tilde{y}^j}\kappa^j(\eta)\dd\eta
\end{equation}

\noindent holds, which by the virtue of (\ref{eq:kappa_connection}) implies $\tilde{y}^i<\tilde{y}^j$. By the properties of $\kappa(y)$ and (\ref{eq:kappa_connection}) it follows that the curvatures are related via
\begin{equation}
\label{eq:kappa_connection2}
\kappa^i(\tilde{y}^i)>(1+\varepsilon)\kappa^j(\tilde{y}^j),
\end{equation}

\noindent because $\tilde{y}^i<\tilde{y}^j$. From this observation and the positivity of all the involved quantities we conclude that 
\begin{equation}
\label{eq:kappa_connection3}
\frac{\tilde{y}^j}{(1+\varepsilon)\kappa^j(\tilde{y}^j)}>\frac{\tilde{y}^i}{\kappa^j(\tilde{y}^i)}.
\end{equation}

We switch to the parametrization of $\bar\Gamma$ with respect to the turning angle $\gamma$. Based on (\ref{eq:DE_kappa_gamma}) the chain rule yields, that the $\bar{A}$ area under $\bar\Gamma$ can be computed as
\begin{equation}
\label{eq:area}
\bar{A}=\int_P^Qy\cos\gamma\dd s=\int_{\frac{\pi}{2}}^0\frac{y}{\kappa}\cos\gamma\dd \gamma.
\end{equation}

As we have demonstrated in (\ref{eq:kappa_connection3}), the argument of the integral in the r.h.s of (\ref{eq:area}) is smaller for $\bar\Gamma_i^*$ than for $\bar\Gamma_j^*$, and this holds for any $\gamma\in(0,\pi/2)$, whence we conclude 
\begin{equation}
\frac{1}{1+\varepsilon}\bar{A}^j=\int_{\frac{\pi}{2}}^0\frac{y^j}{(1+\varepsilon)\kappa^j}\cos\gamma\dd \gamma>\int_{\frac{\pi}{2}}^0\frac{y^i}{\kappa^i}\cos\gamma\dd \gamma=\bar{A}^i.
\end{equation}

Finally we apply (\ref{eq:c_connection}) to obtain
\begin{equation}
\frac{\bar{A}^j}{\hat{c}_1^j}>\frac{\bar{A}^i}{\hat{c}_1^i}.
\end{equation}

As a steady state $\Gamma^*$ curve possesses $D_2$ symmetry $A=4\bar{A}$ follows, so we are left with the conclusion that 
\begin{equation}
\frac{\hat{c}_1^i}{A^i}>\frac{\hat{c}_1^j}{A^j},
\end{equation}

\noindent which is exactly the monotonicity of the $c_1(\hat{c}_1)$ function. This proves that $M$ is injective, as different elements in $\chi_q$ cannot be mapped to identical values. It is also worthy to note, that for all $\hat{c}_1\in\chi_q$ the area is obviously positive thus $c_1(\hat{c}_1)$ is a positive, monotonic, continuous function.
\end{proof}

\begin{Lemma}
 Map M is surjective.
\end{Lemma}
\begin{proof}
To prove surjectivity we have to investigate the limits of $c_1(\hat{c}_1)$ as $\hat{c}_1$ is varied. First we turn to investigate the limit as $\hat{c}_1\rightarrow 0$ ($q$ is still fixed). From Section 3 we know, that the curvature along $\bar\Gamma$ is maximal at point P ($\kappa(0)$) with $\kappa(0)=\hat{c}_1^{-1}$ and it is minimal at point Q with $\kappa(\bar{y})=\hat{c}_1^{-1}(1-2q\bar{y}\phi(\bar y))=\hat{c}_1^{-1}(1-2q\bar{y}\hat{c}_1)$. Curvature of any planar curve is the reciprocal of the $r$ radius of its osculating circle. It provides an estimate on the area of the curve via $
r_{\min}^2\pi<A<r_{\max}^2\pi$, where $r_{\min}$ and $r_{\max}$ are the minimal and maximal radii of the osculating circles along the curve, respectively. Putting it together and we obtain the following inequality
\begin{equation}
\label{ineq01}
\cfrac{\hat{c}_1}{\pi}\left(\cfrac{\hat{c}_1}{1-2q\bar{y}\hat{c}_1}\right)^{-2}<\frac{\hat{c}_1}{A(\hat{c}_1)}<\frac{\hat{c}_1}{\pi}\hat{c}_1^{-2}
\end{equation}

\noindent Recall that $\bar{y}\leq y_0$, hence Lemma 3.1 yield, that at a fixed $q$ the value of $\bar{y}$ is finite. It means that both the lower and the upper expression in the above inequality approach $+\infty$ as $\hat{c}_1\rightarrow 0$. We conclude 
\begin{equation}
\label{eq:limits_0}
\lim_{\hat{c}_1\rightarrow 0}\frac{\hat{c}_1}{A(\hat{c}_1)}=+\infty.
\end{equation}

Finally we investigate the $\hat{c}_1\rightarrow \hat{c}_{\crit}$ limit. As $\hat{c}_{\crit}$ is finite it is enough to investigate the $\bar{A}(\hat{c}_1)$ area in the limit. We consider the already used identity between the curvature and and arch length. Taking again the parametrization with respect to $\gamma$ we write
\begin{equation}
\label{eq:kappa_gamma}
\kappa(\gamma)=- \left(\frac{\dd S(\gamma)}{\dd\gamma}\right)^{-1},
\end{equation}

\noindent where $S(\gamma)$ is the arch length between point P and the point with turning angle $\gamma$. As at $\hat{c}_1=\hat{c}_{\crit}$ the curvature at point Q vanishes we conclude, that
\begin{equation}
\label{eq:limit_S}
\lim_{\gamma\rightarrow 0}\frac{\dd S(\gamma)}{\dd\gamma}=\lim_{\gamma\rightarrow 0}\frac{1}{\kappa(\gamma)}=\infty.
\end{equation}

Thus the curve is unbounded. As the area $\bar{A}$ under $\bar\Gamma$ can be computed from the arc length ($y$ is finite!) we obtain
\begin{equation}
\label{eq:limits_inf}
\lim_{\hat{c}_1\rightarrow\hat{c}_{1,\crit}}S=\lim_{\hat{c}_1\rightarrow\hat{c}_{1,\crit}}A=\infty,
\end{equation}

\noindent which provides the required limit as
\begin{equation}
\label{eq:limits_crit}
\lim_{\hat{c}_1\rightarrow \hat{c}_{\crit}}\frac{\hat{c}_1}{A(\hat{c}_1)}=0.
\end{equation}

It means, the range of $M$ is indeed $\xR^+$ and based on the injectivity part of the proof the preimage is precisely $\chi_q$. 
\end{proof}

\begin{Rem}
The arguments above show, that the solution curve $\Gamma$ is compact for any $0<\hat{c}_1<\hat{c}_{1,\crit}$ and $0\leq\hat c_2<\infty$. Using (\ref{eq:solvability_3}) we see that compactness holds for parameters $\sqrt{\hat c_1\hat c_2}<\Psi$. $\Gamma$ is not compact iff $\sqrt{\hat c_1\hat c_2}=\Psi$
\end{Rem}

\begin{proof}[Proof of Theorem \ref{thm2}]
As $M$ is injective and surjective we conclude that it must be \emph{one-to-one and onto}. This means,that the nonlocal equation in (\ref{eq:steady_state}) produces a unique, compact solution among smooth curves for any positive $c_1$ and $c_2$. 
\end{proof}

\begin{Rem}
We proved that time-invariant, smooth curves under the flow in (\ref{eq:geom_PDE}) are uniquely determined by the $c_1$ and $c_2$ parameters and they possess $D_2$ symmetry. Observe that uniqueness stems from the linear ODE obtained for $\cos\gamma(y)$, which ensures uniqueness for the $\kappa(y)$ curvature function. It seems that other laws either to the abrasion or to the friction term might destroy this linearity. It seems even more probable, that other non-local quantities (instead of $A$ in (\ref{eq:geom_PDE})) would either lead to collapse of the injectivity or the surjectivity of the map $M$. 
\end{Rem}

\section{Conclusion and outlook}
\label{sec:conclusions}
Based on physical intuition a model of ooid-growth in 2D was introduced in \cite{Sipos2018}. There a remarkable similarity between model predictions and natural shapes were found. Here we rigorously prove existence and uniqueness of time-invariant shapes under the flow. Investigation of a broad class of related flows is an interesting future project, as well as the investigation of the spatial version of the model. The practical significance of the presented results lies in the unique relation between the physical relation of the time-invariant shape and the model parameters. It implies that pure observation of ooid shapes and cross sections can be directly used to deduce features of the physical environment that formed the particle. Hence, this work motivates deeper understating of the connections between the model parameters and physical characteristics.

\section*{Acknowledgment}
I am indebted to the anonymous referee for his/her comments which significantly improved the manuscript. I thank G\'abor Domokos for his idea to investigate the model presented in the paper and the fruitful discussions about ooids. Support of the NKFIH grant K 119245 and grant BME FIKP-V\'IZ by EMMI is gratefully acknowledged.

\end{document}